\documentclass{article}
\usepackage{a4wide}
\usepackage[utf8]{inputenc}
\usepackage{amsmath}
\usepackage{amsthm}
\usepackage{etoolbox}
\usepackage[usenames]{color}
\usepackage[usenames]{xcolor}
\usepackage{amsfonts}
\usepackage{graphicx}      
\usepackage{tikz-cd}
\usetikzlibrary{babel}
\usepackage{natbib}
\usepackage{comment}
\usepackage{url}

\patchcmd{\subequations}{}%
{}{}{}

\newcommand{\abs}[1]{\left\lvert #1\right\rvert}

\newcommand{\ba}{\mathsf{a}}
\newcommand{\bb}{\mathsf{b}}

\newcommand{\calL}{\mathcal{L}}

\newcommand{\Vn}{\mathbf{n}}

\newcommand{\VE}{\mathbf{E}}
\newcommand{\VB}{\mathbf{B}}

\newcommand{\R}{\mathbb{R}}

\newcommand{\N}{\mathbb{N}}

\newcommand{\norm}[2]{\left\lVert #1\right\rVert_{#2}}
\definecolor{orange}{rgb}{1,0.4,0}
\definecolor{green}{rgb}{0,0.4,0}

\newcommand{\curl}{\operatorname{\boldsymbol{\mathrm{curl}}}}

\newcommand{\Div}{\operatorname{\mathrm{div}}}
\newcommand{\grad}{\nabla}

\usepackage{stmaryrd}
\usepackage{mathbbol}
\usepackage{amssymb}             

\DeclareSymbolFontAlphabet{\amsmathbb}{AMSb}%

\newcommand{\dd}{\mathrm{d}}

\newcommand{\ddb}[1]{\mathrm{d}_{\Gamma}}

\newcommand{\Tr}{\mathrm{Tr}} 

\newcommand{\HL}[1]{\Delta_{\mathrm{HL},#1}}

\newcommand{\XcZ}[1]{X^{#1}_{0;,0}}
\newcommand{\XZc}[1]{X^{#1}_{0;0,}}
\newcommand{\XX}[1]{X^{#1}}
\makeatletter
\newsavebox{\@brx}
\newcommand{\llangle}[1][]{\savebox{\@brx}{\(\m@th{#1\langle}\)}%
  \mathopen{\copy\@brx\mkern2mu\kern-0.9\wd\@brx\usebox{\@brx}}}
\newcommand{\rrangle}[1][]{\savebox{\@brx}{\(\m@th{#1\rangle}\)}%
  \mathclose{\copy\@brx\mkern2mu\kern-0.9\wd\@brx\usebox{\@brx}}}
\makeatother

\renewcommand{\R}{\amsmathbb{R}}

\newcounter{counter}
\numberwithin{counter}{section}
\newtheorem{theorem}[counter]{Theorem}

\newtheorem{corollary}[counter]{Corollary}
\newtheorem{lemma}[counter]{Lemma}

\newtheorem{rmk}{Remark}

\title{A Unified Space-Time Variational Framework for Wave Equations on
Differential Forms}
\author{C.~Urz\'ua-Torres. A.~Wisse\footnote{Delft Institute of Applied
Mathematics, TU Delft ({\tt C.A.UrzuaTorres@tudelft.nl},
{\tt A.C.Wisse@tudelft.nl}).
This work was funded by the Dutch Research Council (NWO) under the NWO-Talent
Programme Veni with the project number VI.Veni.212.253.}}

\begin{document}
\maketitle

\begin{abstract}
We revisit a standard space--time variational formulation of the scalar wave
equation in the setting of differential forms. For this, we use as spatial
operator the nonnegative, self-adjoint (relative) \emph{Hodge--Laplacian} on a
bounded Lipschitz domain. Its spectral decomposition reduces the wave problem in
any form degree to the same family of scalar ordinary differential equations that
underlies existing space-time analysis. We illustrate how this can be used to
transfer results from the scalar wave equation with Dirichlet boundary conditions
to other wave-type problems. In particular, we connect the framework to perfectly
electrically conducting Maxwell problems.
\end{abstract}

\section{Introduction}
\label{sec:introduction}
The present work is motivated by the development of stable space-time Galerkin
methods for wave-type problems. Classical continuous theory for wave equations
typically relies on a separation of variables or treats time strictly as an
evolution parameter, which naturally leads to semi-discrete numerical schemes
such as the method of lines. While this approach is well established (see, e.g.,
\cite{Hug12, Bak76}), it is not possible to straightforwardly use this in
combination with parallelization and adaptivity in space and time simultaneously.
This is the promise of space-time methods, which also naturally accommodate
moving domains.

Wave equations naturally arise when modeling a broad range of physical
processes, including acoustics, electromagnetism, and the dynamics of
elastic materials; see for example \cite{Gri17, Par03}. The
associated spatial differential operators and boundary conditions depend on
the application. Nevertheless, a common analytical mechanism emerges whenever
the spatial operator admits a nonnegative self-adjoint realization with
discrete spectrum: the expansion in spatial eigenfunctions reduces the wave
equation to scalar problems
\begin{equation}\label{eq:FourierIntro}
 a_m''(t)+\mu_m a_m(t)=f_m(t).
\end{equation}
This is the strategy used in Steinbach and Zank's analysis of the scalar
wave equation \cite{StZ20}, and many other works in space-time formulations
afterwards, see for example \cite{FLS26, HLU25, LSZ22, STZ26, StZ22}. The purpose
of this paper is to use differential forms to make the transfer of these results
explicit for more general wave equations in a unified space-time setting.

Differential forms provide a natural language that encompasses scalar waves and
vector systems and makes the role of boundary conditions explicit. The Hodge
wave equation and related numerical methods have been studied in
\cite[Sec.~1]{WuB21}. We also refer to \cite[Sec.~3 and 6.2]{AFW10} for the
associated Hilbert-complex framework and boundary conditions. However, although
the key ingredients are known and the problem has been considered numerically,
to the best of the author's knowledge, this is the first time it has been
studied in the context of space-time Galerkin methods.

Let us begin by introducing our main character: the Hodge-Laplacian
\begin{equation*}
  \HL{k}:= \delta \dd + \dd \delta
\end{equation*}
over $k$-forms, with $0\leq k \leq n$. Here, $\dd$ denotes the exterior
derivative and $\delta$ is the coderivative.

Let $\Omega\subset \R^n$, $n = 1,2,3$, be a bounded contractible Lipschitz
domain. Let $0<T<\infty$, and define the space-time cylinder $Q:=\Omega
\times(0,T)$, with lateral boundary $\Sigma:= \partial\Omega \times[0,T]$.
For sufficiently smooth $A$, we consider the generalized wave equation
\begin{align}\label{eq:general_wave}
  \partial_{tt}A + \HL{k}A &= j, \quad \text{ in } Q\\
  A(x,0) = \partial_t A(x,0) &= 0, \quad \text{ for } x\in\Omega\\
  \Tr(A)  = \Tr(\delta A)&= 0,\quad \text{ on } \Sigma,
\end{align}
where $\Tr$ is the trace operator mapping to the boundary.

When $n=3$ and $k=1$, the differential expression becomes
\begin{equation*}
 \delta\dd A+\dd\delta A
 \quad\longleftrightarrow\quad
 \curl\curl\boldsymbol A-\grad\Div\boldsymbol A=-\Delta\boldsymbol A.
\end{equation*}

Indeed, in three spatial dimensions, our approach simultaneously covers the
scalar Dirichlet problem ($k=0$), the electric-field and magnetic-field wave
problems ($k=1$ and $k=2$), and the scalar Neumann problem ($k=3$). It is worth
noting that this last case is not just a formal extension of the Dirichlet case.
The main difference is that the Neumann zero mode must be retained in the
spectral analysis, reason why our proofs will be slightly different than a
mere extension of those in \cite{StZ20}. For the scalar Neumann problem, an
extended space-time theory was developed in
\cite[Sec.~2.2 and Sec.~4.3, Lemmas~4.1--4.2]{SUT22}. We use
that work as motivation to see how the current unified framework is
useful in other contexts.

\paragraph{Connection with Maxwell's equations.}
In a homogeneous isotropic medium with constant permittivity $\varepsilon>0$
and permeability $\mu>0$, we identify the electric field with a $1$-form $E$;
the magnetic induction with a $2$-form $B$; and the current density with a
$1$-form $J$. With $c^2=(\varepsilon\mu)^{-1}$ and assuming sufficiently
regular data,  Maxwell's equations imply that
\begin{equation}\label{eq:maxwellWaves}
 \begin{aligned}
 E_{tt}+c^2(\delta\dd+\dd\delta)E
   &=-\varepsilon^{-1}(\partial_tJ+c^2\dd\rho),\\
 B_{tt}+c^2(\delta\dd+\dd\delta)B
   &=\varepsilon^{-1}\dd J.
 \end{aligned}
\end{equation}
For a perfectly electrically conducting (PEC) boundary, the conditions
$\Vn\times\VE=0$ and $\Vn\cdot\VB=0$ are precisely $\Tr E=\Tr B=0$. When
the trace operators are well-defined, Maxwell's divergence and evolution
equations further give that
\begin{equation*}
 \Tr(\delta E)=-\varepsilon^{-1}\Tr\rho,
 \qquad \Tr(\delta B)=\mu\Tr J.
\end{equation*}
In particular, homogeneous boundary conditions arise for PEC problems whose
prescribed volume sources vanish near the boundary. For sources reaching the
boundary, the second boundary condition is generally inhomogeneous. Therefore,
the corresponding treatment requires homogenization of boundary data via a
lifting under appropriate regularity and compatibility assumptions. Moreover,
the second-order equations \eqref{eq:maxwellWaves} represent Maxwell's
evolution only for compatible first-order constraints and initial data.
Alternatively, electromagnetic potentials in Lorenz gauge satisfy wave
equations of Hodge type, subject to gauge and source constraints. While we
do not pursue these equivalences in this paper, they motivate why our model
problem is still meaningful for Maxwell's equations.

\section{Notation}
\subsection{Spaces over the domain $\Omega$}
As before, let $\Omega \subset \R^n$, $n=1,2,3$, be a bounded contractible
Lipschitz domain. We follow the notation from \cite{AFW06}. For $0\leq k\leq n$,
we define $\Lambda^k(\Omega)$ as the space of all smooth differential $k$-forms
on $\Omega$.
We define the $L^2$-inner product of two $k$-forms as
\begin{equation*}
(\omega,\eta)_{L^2\Lambda^k(\Omega)} := \int_{\Omega} \omega \wedge \star \eta,
\end{equation*}
where $\wedge$ denotes the so-called exterior product, as defined in
\cite[Equation (2.1)]{AFW06}.
The Hilbert space $L^2\Lambda^k(\Omega)$ is defined as the completion of
$\Lambda^k(\Omega)$ in the corresponding $L^2\Lambda^k(\Omega)$-norm.
Next, we introduce the spaces
\begin{align*}
H \Lambda^k(\Omega) :=\{\omega \in L^2 \Lambda^k(\Omega) |\; \dd \omega \in
L^2 \Lambda^{k+1}(\Omega)\},
\end{align*}
with norm $\|\omega\|^2_{H \Lambda^k(\Omega)} := \|\omega\|^2_{L^2 \Lambda^k
(\Omega)} + \|\dd \omega\|^2_{L^2 \Lambda^{k+1}(\Omega)}$; and
\begin{align*}
H^* \Lambda^k(\Omega) &:=\{\omega \in L^2 \Lambda^k(\Omega) |\; \delta \omega
\in L^2 \Lambda^{k-1}\},
\end{align*}
with norm $\|\omega\|^2_{H^* \Lambda^k(\Omega)} := \|\omega\|^2_{L^2 \Lambda^k
(\Omega)} + \|\delta \omega\|^2_{L^2 \Lambda^{k-1}(\Omega)}$.

We also consider the subspaces
\begin{align*}
\mathring{H} \Lambda^k(\Omega)&:=\{\omega \in H \Lambda^k(\Omega)|\;
\Tr(\omega) = 0\},\\
V^k(\Omega) &:= H \Lambda^k(\Omega)\cap H^* \Lambda^k(\Omega),\\
\mathring{V}^k(\Omega)&:= \mathring{H}\Lambda^k(\Omega)\cap H^*\Lambda^k(\Omega).
\end{align*}
Let us define the bilinear form
\begin{align}\label{eq:bbk}
 \bb_k(\omega,\beta):= (\dd\omega,\dd\beta)_{L^2\Lambda^{k+1}(\Omega)}
 + (\delta\omega,\delta\beta)_{L^2\Lambda^{k-1}(\Omega)}, \qquad
 \text{ for }\omega,\beta\in V^k(\Omega).
\end{align}
We equip $V^k(\Omega)$ with the following inner-product
\begin{equation*}
(\omega,\beta)_{V^k(\Omega)}:= (\omega,\beta)_{L^2 \Lambda^k(\Omega)}+
 \bb_k(\omega,\beta),\qquad \text{ for }\omega,\beta\in V^k(\Omega),
\end{equation*}
which induces the semi-norm $|\omega|^2_{V^k(\Omega)}:=  \|\dd \omega\|^2_{L^2
\Lambda^{k+1}(\Omega)} +  \|\delta \omega\|^2_{L^2 \Lambda^{k-1}(\Omega)}$
and norm $$\|\omega\|^2_{V^k(\Omega)} := |\omega|^2_{V^k(\Omega)} + \|\omega
\|^2_{L^2\Lambda^k(\Omega)}.$$
\begin{rmk}\label{remark:equiv_norm_seminorm}
Note that due to the Poincaré inequality \cite[Theorem 4.6]{Arn18},
the $V^k(\Omega)$-norm and semi-norm are equivalent on $\mathring{V}^k(\Omega)$
when $0\leq k<n$. Similarly, we can leave out the first term in the
inner-product for $0\leq k<n$.
\end{rmk}

\subsection{Spaces over the time interval $(0,T)$}
We now define the spaces over the temporal domain, needed to hereafter build
spaces over the space-time cylinder $Q$. We consider the following spaces with
initial/terminal conditions
\begin{align*}
  H^1_{0,}(0,T) &:= \{v\in H^1(0,T): v(0) = 0\},\\
  H^1_{,0}(0,T) &:= \{v\in H^1(0,T): v(T) = 0\}.
\end{align*}
Here, $H^1(0,T)$ is the standard first-order Sobolev space. These spaces are
endowed with the inner product
\begin{equation*}
(u,v)_{H^1(0,T)}:= \int_0^T \partial_t u(t) \partial_t v(t) dt,
\end{equation*}
and induced semi-norm
\begin{equation*}
|u|^2_{H^1(0,T)}:= \|\partial_t u\|^2_{L^2(0,T)} = \int_0^T |\partial_t u(t)|^2
dt.
\end{equation*}

\subsection{Spaces over the space-time cylinder $Q$}
We define the following spaces:
\begin{align}
\XX{k}&:= H^1(0,T;L^2 \Lambda^k(\Omega))\cap L^2(0,T; V^k(\Omega)),
\label{eq:H}\\
\XZc{k}&:= H^1_{0,}(0,T;L^2 \Lambda^k(\Omega))\cap
L^2(0,T; \mathring{V}^k(\Omega)),\label{eq:HcZ}\\
\XcZ{k}&:= H^1_{,0}(0,T;L^2 \Lambda^k(\Omega))\cap L^2(0,T; \mathring{V}^k(
\Omega)),\label{eq:HZc}
\end{align}
with semi-norm
\begin{align*}
  |\omega|^2_{\XX{k}} := &\int_0^T \lbrace(\partial_t\omega,\partial_t\omega)_{
  L^2\Lambda^k(\Omega)} + \bb_k(\omega,\omega)\rbrace dt
\end{align*}
and norm $
\|\omega\|^2_{\XX{k}} := |\omega|^2_{\XX{k}}+\int_0^T (\omega,\omega)_{L^2
\Lambda^k(\Omega)}dt.$

For spaces of this type, we have the following results.
\begin{lemma}\label{lem:isom}
The following isometries hold:
\begin{itemize}
\item $H^1(0,T; L^2\Lambda^k(\Omega))$ is isometric to the Hilbert tensor
  product $H^1(0,T)\hat{\otimes}L^2\Lambda^k(\Omega)$,
\item $H_{a}^1(0,T; L^2\Lambda^k(\Omega))$ is isometric to the Hilbert tensor
product $H_{a}^1(0,T)\hat{\otimes} L^2\Lambda^k(\Omega)$, where $a$ is either
replaced by $,0$ or $0,$,
\item Let $F$ be either $V^k(\Omega)$ or $\mathring{V}^k(\Omega)$. Then
$L^2(0,T; F)$ is isometric to the Hilbert tensor product $L^2(0,T)\hat{\otimes}
F$.
\end{itemize}
\end{lemma}
\begin{proof}
The first statement follows from \cite[12.7 Thm 1]{Aub79} and the third
statement follows from \cite[Sect.~12.6 Thm 1]{Aub79}. We get the second
statement by adapting the proof of \cite[Sect.~12.7 Thm 1]{Aub79} so that we
take $H^1_a(0,T)$ as domain of the operator $\partial_t$. The rest of the
proof is identical.
\end{proof}

Note that by a temporal Poincar\'e inequality, we have that the semi-norm $|
\cdot|_{\XX{k}}$ and norm $\|\cdot\|_{\XX{k}}$ are equivalent in the spaces
$\XZc{k}$ and $\XcZ{k}$ ($0\leq k \leq n$).

\section{Spectral Properties of the Hodge-Laplacian}

In this section we present some key spectral properties of the Hodge-Laplacian
that are needed to generalize the results from the Dirichlet scalar wave equation
to the current setting.
In particular, we are interested in the nonnegative self-adjoint operator
$\calL_k:=\Delta_{HL,k}$ associated with the densely defined closed form $\bb_k$
on $\mathring{V}^k$, as introduced in \eqref{eq:bbk}. Hence, its operator domain
is given by
\begin{equation}\label{eq:domainLk}
 D(\calL_k):=\bigl\{u\in \mathring{V}^k:\ \exists f\in L^2\Lambda^k(\Omega)
 \text{ such that }\bb_k(u,v)=(f,v)_{L^2\Lambda^k(\Omega)}\ \forall v\in
 \mathring{V}^k\bigr\}.
\end{equation}
For smooth forms on a smooth $\Omega$, this domain corresponds to the boundary
conditions $\Tr u=0$ and $\Tr(\delta u)=0$. On Lipschitz domains $\Omega$,
\eqref{eq:domainLk} considers the definitions
\cite[eqs.~(1.1)--(1.5)]{MTV05}.

We are now in the position to state the first required spectral result:
\begin{theorem}\label{thm:spectralHL}
Let $\Omega\subset\mathbb R^n$ be a bounded Lipschitz domain and $0\leq k\leq n$.
The operator $\calL_k$ in \eqref{eq:domainLk} is nonnegative, self-adjoint, and has
compact resolvent. Consequently there is an $L^2\Lambda^k(\Omega)$-orthonormal
eigenbasis $\{\phi_m\}_{m\geq1}$, with
\begin{equation}
  \calL_k\phi_m=\mu_m\phi_m, \qquad 0\leq \mu_1 \leq \mu_2 \leq ...
\end{equation}
where each eigenvalue has finite multiplicity, and the eigenvalues accumulate
only at infinity.
\end{theorem}
\begin{proof}
The spectral result follows from the standard spectral theorem for compact
self-adjoint operators, as in \cite[Theorem~1.30]{Ros97}. We
therefore only need to establish the corresponding properties of the
(relative) Hodge--Laplacian on a bounded Lipschitz domain.

First, we recall that the Hodge--Laplacian $\calL_k$ is a nonnegative
self-adjoint operator on $L^2\Lambda^k(\Omega)$, as established in
\cite[eqs.~(1.1)--(1.5)]{MTV05}. Second, we have that the embedding
\begin{equation}
\mathring{V}^k\hookrightarrow L^2\Lambda^k(\Omega)
\label{eq:compact_embedding}
\end{equation}
is compact \cite[Sec.~6.2]{AFW10}.

Now, by the Lax-Milgram theorem applied to the bilinear form
\begin{equation*}
(u,v)_{L^2\Lambda^k(\Omega)}+\bb_k(u,v)
\end{equation*}
on $\mathring V^k$, we have that the resolvent $(I+\calL_k)^{-1}$ is well-
defined and bounded as a mapping from $L^2\Lambda^k(\Omega)$ into $\mathring
V^k$. This and the compactness of \eqref{eq:compact_embedding} imply that
\begin{equation*}
K_k:=(I+\calL_k)^{-1}:L^2\Lambda^k(\Omega)\to \mathring{V}^k
\hookrightarrow L^2\Lambda^k(\Omega)
\end{equation*}
is compact.

Moreover, the fact that $\calL_k$ is nonnegative and self-adjoint, gives that
$K_k$ is self-adjoint, positive, and injective. Therefore, the spectral theorem
for compact self-adjoint operators provides an orthonormal basis $\{\phi_m\}_{m
\geq1}$ of $L^2\Lambda^k(\Omega)$ satisfying
\begin{equation*}
K_k\phi_m=\nu_m\phi_m, \qquad 0<\nu_m\leq1, \text{ for } m\in\N,
\end{equation*}
and such that each eigenvalue has finite multiplicity and $\nu_m\to0$ when $m
\to\infty$.

Consequently, by definition of $K_k$, we have that
\begin{equation*}
\calL_k\phi_m=\mu_m\phi_m,
\qquad
\mu_m=\frac{1}{\nu_m}-1\geq0,
\end{equation*}
and $\mu_m\to\infty$ as $m\to\infty$.
\end{proof}

Along the lines of \cite[Lemma~3.3.5]{Jos17}, one can show:
\begin{lemma}\label{lem:harmonic}
For each $u\in D(\calL_k)$,
\begin{equation*}
 (\calL_ku,u)_{L^2\Lambda^k(\Omega)}=\norm{\dd u}{L^2\Lambda^{k+1}(\Omega)}^2
 +\norm{\delta u}{L^2\Lambda^{k-1}(\Omega)}^2\geq0,
\end{equation*}
and $\calL_k u=0$ if and only if $\dd u=\delta u=0$.
\end{lemma}
\begin{proof}
By \eqref{eq:domainLk}, $(\calL_k u,v)_{L^2\Lambda^k(\Omega)}=(f,v)_{L^2\Lambda^k
(\Omega)}=\bb_k(u,v)$ for all $v\in \mathring{V}^k$.

Moreover, if $\dd u=\delta u=0$, then $\bb_k(u,v)=0$ for all $v\in\mathring{V
}^k$. Therefore $\calL_k u=0$. The converse follows from the identity.
\end{proof}

With these ingredients, we can characterize the dimension of the kernel of
the Hodge-Laplacian.
\begin{corollary}\label{cor:kernelLk}
Let $\Omega$ be a bounded Lipschitz domain. If, in addition, $\Omega$ is
contractible, then
\begin{equation}\label{eq:kernel}
 \dim\ker \calL_k=\begin{cases}0,&0\leq k<n,\\1,&k=n.\end{cases}
\end{equation}
In particular, all eigenvalues are strictly positive for $k<n$. For $k=n$
exactly one eigenvalue (counting multiplicity) is zero, and all other
eigenvalues are strictly positive.
\end{corollary}
\begin{proof}
By Lemma~\ref{lem:harmonic}, the kernel of the Hodge-Laplacian $\calL_k$ is
\begin{equation}
\ker \calL_k =\left\{u\in \mathring{V}^k: du=0,\ \delta u=0 \right\}.
\label{eq:kernelHL}
\end{equation}

This is precisely the space of relative harmonic $k$-forms, denoted by
$\mathring{\mathfrak H}^k$ in \cite[Sec.~6.2, Eq.~(80)]{AFW10}.
From there, we also know that
\begin{equation*}
\dim\mathring{\mathfrak H}^k =b_{n-k}(\Omega),
\end{equation*}
where $b_j(\Omega)$ denotes the $j$-th Betti number of $\Omega$, i.e., the
dimension of its $j$-th de Rham cohomology space. Now, because these spaces are
trivial for a contractible domain $\Omega$ except when $j=0$, where they consist
of the constant functions \cite[p.~18]{AFW06}, we obtain that
\begin{equation}
b_j(\Omega)=
\begin{cases}
1, & j=0,\\
0, & 1\leq j\leq n.
\end{cases}
\label{eq:contractible_betti}
\end{equation}

Combining \eqref{eq:kernelHL}-\eqref{eq:contractible_betti},
we get
\begin{equation}
\dim\ker \calL_k=
\begin{cases}
0, & 0\leq k<n,\\
1, & k=n.
\end{cases}
\label{eq:dimKerLk}
\end{equation}

Since $\calL_k$ is nonnegative by Theorem~\ref{thm:spectralHL},
\eqref{eq:dimKerLk} implies that all its eigenvalues are
strictly positive for $0\leq k<n$. For $k=n$, this tells us that zero is an
eigenvalue of multiplicity one, and all remaining eigenvalues are strictly
positive.

The finite multiplicity of the eigenvalues and their accumulation only at
infinity also follow from Theorem~\ref{thm:spectralHL}.
\end{proof}

It is worth noticing that Corollary~\ref{cor:kernelLk} establishes the reason
why our approach needs to slightly depart from Steinbach and Zank's strategy:
We have to deal with the possibility of a zero eigenvalue. This will also
motivate that the next result considers the full graph inner product.
\begin{lemma}\label{lem:basis_0Vk}
Let $\Omega$ be a bounded Lipschitz contractible domain. Then, the eigenforms
$\{\phi_m\}$ of Theorem~\ref{thm:spectralHL} form an orthogonal basis of
$\mathring{V}^k$, equipped with the full graph inner product.
\end{lemma}
\begin{proof}
Let $\{\phi_m\}_{m\geq1}$ be the basis of $L^2\Lambda^k(\Omega)$ from
Theorem~\ref{thm:spectralHL}, with corresponding eigenvalues $\mu_m$, i.e.,
$L^2 \Lambda^k(\Omega) = \overline{span\{\phi_m\}}^{\|\cdot\|_{L^2\Lambda^k(
\Omega)}}$. Writing $W: = \overline{span \{\phi_m\}}^{\|\cdot\|_{V^k(\Omega)}}$,
we want to show that $W=\mathring{V}^k(\Omega)$.

Since $\mathring{V}^k(\Omega)$ is a Hilbert space, and $W$ is a closed subset
of $\mathring{V}^k(\Omega)$ by definition, it follows from
\cite[Theorem 3.13]{Nee22} that
\begin{equation*}
  \mathring{V}^k(\Omega) = W \oplus W^{\perp}.
\end{equation*}
Therefore, if we can show that $W^{\perp} = \{0\}$, we are done.

Let $f\in W^{\perp}$. Then we have that
\begin{equation*}
(f, \phi_m )_{V^k(\Omega)} = (f, \phi_m )_{L^2\Lambda^k(\Omega)}
+ \bb_k(f,\phi_m ) = 0,\quad\forall m\geq1,
\end{equation*}
with $\bb_k$ as defined in \eqref{eq:bbk}.

Since $\phi_m$ are eigenforms of the Hodge-Laplacian $\calL_k$, it follows that
\begin{equation}\label{eq:newInnerProd}
(f, \phi_m )_{V^k(\Omega)} = (f, \phi_m )_{L^2\Lambda^k(\Omega)}
+ \mu_m (f,\phi_m )_{L^2\Lambda^k(\Omega)} = 0,\quad\forall m\geq1.
\end{equation}
Next, given that $\mu_m\geq0$, \eqref{eq:newInnerProd} implies that
\begin{equation*}
(f,\phi_m)_{H^k}=0, \qquad \forall m\geq1.
\end{equation*}
Hence, since $\lbrace\phi_m\rbrace_{m\geq1}$ form a basis of $L^2\Lambda^k(
\Omega)$, it follows that $f = 0$. Therefore, $W^\perp=\{0\}$ and $W=\mathring
V^k$.

Finally, for any $i,m\geq1$, by using \eqref{eq:newInnerProd} and $L^2\Lambda^k
(\Omega)$-orthogonalty, we get
\begin{equation*}
(\phi_i,\phi_m)_{V^k} =(1+\mu_m)(\phi_i,\phi_m)_{L^2\Lambda^k(\Omega)}
=(1+\mu_m)\delta_{im},
\end{equation*}
where $\delta_{im}$ is the Kronecker delta. Thus, the eigenforms are also orthogonal
in $\mathring V^k$ equipped with the full graph inner product.
\end{proof}

\begin{rmk}
When $k<n$, it suffices to use the partial graph inner product given by $\bb_k$
in Lemma~\ref{lem:basis_0Vk}, ·since in that case $\mu_m>0$ for all $m\geq1$.
\end{rmk}

\section{The Variational Formulation and The Transfer of Space-Time Results}
Let $\omega \in \XZc{k}, \beta \in \XcZ{k}$. We write
\begin{equation*}
(\omega,\beta)_{L^2(0,T;L^2 \Lambda^k(\Omega))}:= \int_0^T (\omega,\beta)_{
  L^2\Lambda^k(\Omega)}dt,
\end{equation*}
and define the bilinear form
\begin{align*}
\ba(\omega,\beta):= &-\int_0^T (\partial_t\omega, \partial_t\beta)_{L^2\Lambda^k(\Omega)}dt
+ \int_0^T \bb_k(\omega,\beta)dt,
\end{align*}
and the functional
\begin{equation*}
\ell(\beta):= (j, \beta)_{L^2(0,T;L^2\Lambda^k(\Omega))}.
\end{equation*}

For $j\in L^2(0,T;L^2\Lambda^k(\Omega))$, The variational formulation of
\eqref{eq:general_wave} is given by: Find $A\in \XZc{k}$ such that
\begin{equation}\label{eq:ST_varWE}
  \ba(A,\beta) = \ell(\beta) \qquad \text{ for all } \beta \in \XcZ{k}.
\end{equation}
We can now illustrate how the developed framework allows us to transfer
some results from the scalar Dirichlet wave problem to the more general
case of \eqref{eq:general_wave}.

\subsection{Existence and Uniqueness}
In the theorem below, we show existence and uniqueness of solutions to
\eqref{eq:ST_varWE} with right-hand side in $L^2$ in space and time. This is
the more general analogue to \cite[Theorem~5.1]{StZ20} and hence also uses
techniques based on \cite{Lad85}.
\begin{theorem}\label{thm:ex_uni_L2}
Let $j\in L^2(0,T;L^2\Lambda^k(\Omega))$ be given. Then there exists a unique
solution $A\in\XZc{k}$ of \eqref{eq:ST_varWE} satisfying
\begin{equation}\label{eq:L^2_estimate}
  |A|_{\XX{k}}\leq \frac{T}{\sqrt{2}}\|j\|_{L^2(0,T; L^2\Lambda^k(\Omega))}.
\end{equation}
\end{theorem}
\begin{proof}
We follow the structure of the proof of \cite[Theorem~5.1]{StZ20}. Let
$\{\phi_m\}_{m\ge1}$ be the orthonormal basis of $L^2\Lambda^k(\Omega)$ consisting
of eigenforms of the Hodge-Laplace operator $\calL_k$, with corresponding
eigenvalues $\{\mu_m\}_{m\geq 1}$, which is also an orthogonal basis of
$\mathring{V}^k(\Omega)$, as given in Theorem \ref{thm:spectralHL} and Lemma
\ref{lem:basis_0Vk}.

Also, let $\{v_l\}_{l\geq 0}$ be the Fourier basis of $H^1_{0,}(0,T)$ from
\cite[Eq. (2.5)]{StZ20}, which is also an orthogonal basis of $L^2(0,T)$.
Then, any $A\in \XZc{k}$ admits the expansion
\begin{equation}\label{eq:ansatz_A}
A(x,t) = \sum_{m=1}^{\infty}\sum_{l=0}^{\infty} a_{m,l} v_l(t) \phi_m(x)
= \sum_{m=1}^{\infty} A_m(t) \phi_m(x).
\end{equation}
We consider \eqref{eq:ansatz_A} to find the solution to \eqref{eq:ST_varWE}.
Hence, we have to find the unknown functions $A_m\in H^1_{0,}(0,T)$. For this,
fix $m\ge1$ and choose test functions of the form
\begin{equation*}
 \beta(x,t) = U(t)\phi_m(x),
\qquad U \in H^1_{,0}(0,T).
\end{equation*}
Inserting this into the variational formulation \eqref{eq:ST_varWE} and using
orthonormality of the eigenforms yields: Find $A_m \in H^1_{0,(0,T)}$ such that
\begin{equation}\label{eq:varODE}
-\int_{0}^T \partial_t A_m(t) \partial_t U(t) dt + \mu_m \int_0^T A_m U(t)
dt = \int_{0}^{T} J_m U(t)dt, \quad \forall U \in H^1_{,0}(0,T),
\end{equation}
where $J_m\in L^2(0,T)$ is such that $J_m(t) = (j,\phi_m)_{L^2\Lambda^k(\Omega)}, t\in(0,T)$,
i.e., they are the coefficients of the Fourier expansion
\begin{equation*}
j(x,t) = \sum_{m=1}^{\infty} \sum_{l=0}^{\infty} j_{m,l} v_l(t) \phi_m(x)
  = \sum_{m=1}^{\infty} J_m(t) \phi_m(x),\quad (x,t)\in Q.
\end{equation*}
From this, we can conclude
\begin{align}
\|j\|^2_{L^2(0,T; L^2\Lambda^k(\Omega))} &= \int_{0}^{T} \|j\|^2_{L^2\Lambda^k(
\Omega)}dt =
\sum_{m=1}^{\infty}\sum_{i=1}^{\infty} (\phi_i, \phi_m)_{L^2\Lambda^k(\Omega)}
\int_{0}^{T} J_i(t) J_m(t)dt\nonumber\\
&= \sum_{m=1}^{\infty} \int_0^T |J_m(t)|^2 dt = \sum_{m=1}^{\infty}\|J_m\|^2_{
L^2(0,T)}.\label{eq:rewrite_j}
\end{align}

Now, solving \eqref{eq:varODE} is equivalent to: For each $m\ge1$, find
\begin{equation}\label{eq:strongODE}
 A_m''+\mu_m A_m=J_m,\qquad A_m(0)=A_m'(0)=0.
\end{equation}
We have that the solution of \eqref{eq:strongODE} is given by
\begin{equation*}
 A_m(t)=\begin{cases}
\dfrac{1}{\sqrt{\mu_m}} \int_0^t\sin(\sqrt{\mu_m}(t-s))J_m(s)\,ds,
     &\mu_m>0,\\
\int_0^t(t-s)J_m(s)\,ds,&\mu_m=0.
 \end{cases}
\end{equation*}
If $\mu_m>0$, differentiation gives $A_m'(t) =\int_0^t\cos(\sqrt{\mu_m}(t-s))
J_m(s)\,ds.$
From this we see that, using Cauchy-Schwarz, we get
\begin{equation}\label{eq:EstimateA_i}
 |A_m'(t)|^2+\mu_m|A_m(t)|^2\leq t\int_0^t|J_m(s)|^2\,ds.
\end{equation}
For $\mu_m=0$, the same estimate follows from $A_m'(t)=\int_0^t J_m(s)\,ds$.

Integrating \eqref{eq:EstimateA_i} and interchanging the order of
integration, we arrive to
\begin{equation}\label{eq:modalestimate}
\int_0^T(|A_m'|^2+\mu_m|A_m|^2)\,dt \leq \frac{1}{2}\int_0^T(T^2-s^2)
|J_m(s)|^2\,ds \leq\frac{T^2}{2}\norm{J_m}{L^2(0,T)}^2.
\end{equation}

Now to show the estimate \eqref{eq:L^2_estimate}, we first
note that
\begin{align*}
|A|^2_{\XX{k}}&= \int_{0}^{T}(\partial_t A,\partial_t A)_{L^2\Lambda^k(\Omega)}
dt + \int_{0}^{T}\bb_k(A,A)dt\\
&= \sum_{i=1}^{\infty}\sum_{m=1}^{\infty}\left[ \int_{0}^{T} \partial_t A_i
\partial_t A_m dt(\phi_i,\phi_m)_{L^2\Lambda^k(\Omega)}+\int_{0}^{T} A_i A_m dt
\bb_k(\phi_i,\phi_m)\right].
\end{align*}
Using \eqref{eq:modalestimate} and orthogonality of the eigenforms, we obtain
\begin{align*}
|A|^2_{\XX{k}} &= \sum_{m=1}^{\infty} \left[\int_{0}^{T} |\partial_t A_m|^2+
\mu_m \int_{0}^T|A_m|^2 dt\right]
\leq \frac{1}{2}T^2 \sum_{m=1}^{\infty}\|J_m\|^2_{L^2(0,T)} =
\frac{1}{2}T^2\|j\|^2_{L^2(0,T;L^2\Lambda^k(\Omega))},
\end{align*}
where the last step follows from \eqref{eq:rewrite_j}.
From this we get, by orthogonality of the eigenforms and the
equivalence of $|\cdot|_{X_k}$ and $\|\cdot\|_{X_k}$ on the Ansatz space, that
indeed the series \eqref{eq:ansatz_A} converges in $X^k_{0;0,}$. Moreover, by
standard arguments, one can show that the resulting limit satisfies
\eqref{eq:ST_varWE}.

For uniqueness, let $A$ solve \eqref{eq:ST_varWE} with $j\equiv0$. Testing
against $U(t)\phi_m$, $U \in H^1_{,0}(0,T)$ implies that each scalar coefficient
satisfies the distributional homogeneous version of \eqref{eq:strongODE}.
Hence $A_m=0$ for all $m\geq1$, and $A=0$.
\end{proof}
Moreover, we have the following result regarding the bilinear form from
\eqref{eq:ST_varWE}.
\begin{lemma}\label{lemma:bddness_generalized}
The bilinear form $\ba(\cdot,\cdot): \XZc{k}\times \XcZ{k} \to \R$ from
\eqref{eq:ST_varWE} is bounded,
i.e. there exists a constant $\alpha_a>0$ such that
\begin{equation*}
  |\ba(u,v)|\leq \alpha_a\|u\|_{\XX{k}}\|v\|_{\XX{k}}
\end{equation*}
for all $u\in \XZc{k}$ and $v\in \XcZ{k}$. 
\end{lemma}
\begin{proof}
  The result follows from the Cauchy-Schwarz inequality.
\end{proof}

\subsection{Limits of inf-sup stability}
After showing existence and uniqueness of solutions to \eqref{eq:ST_varWE} with
right hand side in $L^2(0,T;L^2\Lambda^k(\Omega))$, one can wonder if the space
of admissible right hand side functions can be enlarged. A natural choice from the
space-time finite element perspective would be $j\in [\XcZ{k}]'$. However, it
was shown in \cite[Theorem 1.1]{StZ22} that for the scalar Dirichlet wave problem,
this choice cannot guarantee existence of a unique solution. In the following, we
extend this result and prove that this is the case for any $0\leq k\leq n$.
\begin{theorem}\label{thm:no_uniqueness}
There does not exist a constant $C>0$ such that each right hand side $j\in
[\XcZ{k}]'$ and corresponding solution $A$ to \eqref{eq:ST_varWE} satisfy
\begin{equation}\label{eq:estimate_sol}
  |A|_{\XX{k}}\leq C \|j\|_{[\XcZ{k}]'}.
\end{equation}
In particular, there does not exist a constant $C_S>0$ such that the inf-sup
condition
\begin{equation}\label{eq:inf-sup}
\forall A\in \XZc{k}:\; \sup_{0\neq w \in \XcZ{k}}\frac{|\ba(A,w)|}{\|w\|_{\XX{
k}}}\geq C_S \|A\|_{\XX{k}}
\end{equation}
holds.
\end{theorem}
\begin{proof}
In order to highlight how the machinery extends, we follow the proof from
\cite[Theorem 4.2.23]{Zan19}.
Let $\{\phi_m\}_{m\geq1}$ be the eigenforms from \ref{thm:spectralHL} with
corresponding eigenvalues ${\mu_m, m\geq1}$. For $i\in\N$ and without loss of
generality, we choose the pair $(\mu_i,\phi_i)$ such that $\mu_i>0$ and define
for $x\in\Omega$ and $t\in[0,T]$
\begin{align*}
A_i(x,t) &= \int_0^t s \sin(\sqrt{\mu_i}s)ds \phi_i(x) =: f_i(t) \phi_i(x),\\
j_i(x,t) &= \partial_{tt}A_i + \HL{k} A_i = 2 \sin(\sqrt{\mu_i}t)\phi_i(x).
\end{align*}
We will see that this choice is the counterexample that breaks
\eqref{eq:estimate_sol} and \eqref{eq:inf-sup}.

First, note that $A_i\in \XZc{k}$. Moreover, it follows that
\begin{align*}
|A_i|^2_{\XX{k}} = &\int_0^T (\partial_t A_i, \partial_t A_i)_{L^2\Lambda^k}dt
+ \int_0^T \bb_k(A_i,A_i) dt.
\end{align*}
Using \cite[Eq 2.8]{AFW06} and the fact that we can write $A_i = f_i(t)\phi_i(x
)$, we get that
\begin{align*}
|A_i|^2_{\XX{k}} = \|\partial_t f_i\|^2_{L^2(0,T)}\|\phi_i\|^2_{L^2\Lambda^k(
\Omega)}+ \|f_i\|^2_{L^2(0,T)}(\calL_k \phi_i, \phi_i)_{L^2 \Lambda^k(\Omega)}.
\end{align*}
Now, we will use that the $\phi_i$ are eigenforms of the Hodge-Laplacian
$\calL_k$, and combine this with the fact that the eigenforms form an orthonormal
basis of $L^2 \Lambda^k(\Omega)$, see Theorem \ref{thm:spectralHL}. This yields
that
\begin{align}\label{eq:semi-normEquiv}
|A_i|^2_{\XX{k}} = \|\partial_t f_i\|^2_{L^2(0,T)}+ \mu_i \|f_i\|^2_{L^2(0,T)},
\end{align}
and the same calculations as in the proof of \cite[Thm 4.2.6]{Zan19} give that
\begin{equation*}
|A_i|_{\XX{k}} \to \sqrt{\frac{T^3}{3}}\quad\text{as } i\to\infty.
\end{equation*}
Next, we compute $\|j_i\|_{[\XcZ{k}]'}$. From the definition of the norm,
we have that
\begin{equation*}
\|j_i\|_{[\XcZ{k}]'}:= \sup_{0\neq w \in \XcZ{k}}\frac{\abs{\langle j_i, w
\rangle}}{|w|_{\XX{k}}}.
\end{equation*}
Write $w(x,t) = \sum_{l=1}^{\infty}w_l(t) \phi_l(x)$, where $\{\phi_l\}$
are the basis from Theorem~\ref{thm:spectralHL}. Then,
\begin{align*}
\langle j_i, w\rangle= (j_i,w)_{L^2(0,T;L^2\Lambda^k(\Omega))}
&= \sum_{l=1}^{\infty} \int_0^T w_l(t)2\sin(\sqrt{\mu_i}t)dt
(\phi_l, \phi_i)_{L^2 \Lambda^k(\Omega)}\\
&= \int_0^T w_i(t)2\sin(\sqrt{\mu_i}t)dt\, \|\phi_i\|^2_{L^2 \Lambda^k(\Omega)}.
\end{align*}
By orthonormality of $\phi_l$ and Cauchy-Schwarz, we get
\begin{align}\label{eq:CSestimate}
\abs{\langle j_i, w\rangle}
= \abs{\int_0^T w_i(t)2\sin(\sqrt{\mu_i}t)dt} \leq 2\sqrt{T}\norm{w_i}{L^2(0,T)}.
\end{align}
On the other hand, by the same arguments used to derive \eqref{eq:semi-normEquiv},
we have that
\begin{align}\label{eq:lower_bound}
|w|_{X^k}^2 &=\|\partial_t w\|_{L^2(0,T;L^2\Lambda^k(\Omega))}^2 +\|\dd w\|_{L^2(Q
)}^2 +\|\delta w\|_{L^2(Q)}^2 \nonumber\\
&=\sum_{l=1}^{\infty} \left(\|w_l'\|_{L^2(0,T)}^2 +\mu_l\|w_l\|_{L^2(0,T)}^2
\right)\geq \mu_i\|w_i\|_{L^2(0,T)}^2.
\end{align}
Finally, combining \eqref{eq:lower_bound} with \eqref{eq:CSestimate}, we obtain
\begin{equation*}
\|j_i\|_{[\XcZ{k}]'}= \sup_{0\neq w \in \XcZ{k}}\frac{\langle j_i, w
\rangle}{|w|_{\XX{k}}} \leq \sup_{0\neq w \in \XcZ{k}}\frac{\frac{2\sqrt{T}}{
\sqrt{\mu_i}}|w|_{X^k}.}{|w|_{\XX{k}}} \to 0 \text{ as } i \to \infty.
\end{equation*}

We conclude the proof by using that $|\cdot|_{X_k}$ and $\|\cdot\|_{X_k}$ are
equivalent in the considered spaces, and by noting that, since $\ba(A_i,w)=
\langle j_i,w\rangle$ for every $w\in \XcZ{k}$, these two limits contradict the
existence of a positive uniform inf-sup constant. Hence neither
\eqref{eq:estimate_sol} nor \eqref{eq:inf-sup} can hold.
\end{proof}

It is worth pointing out that in the same way that our framework allowed us to
transfer the results from \cite{Zan19} to our generalized wave equation, it also
paves the way to accommodate for $j\in[\XcZ{k}]'$ by transitioning to a
distributional framework, as demonstrated in \cite[Sec.~3]{StZ22} for the scalar
wave equation.

\section{Discussion and Outlook}
We presented the ingredients to transfer main space-time results for the wave
equation, as initiated in \cite{Zan19, StZ20, StZ22}. Indeed, we described the
precise Hodge-Laplacian operator and energy spaces for which such extension is
valid, including the case $k=n$. The advantage of this unified framework is not
only that it unveils what would be the corresponding electromagnetic counterpart
to the proposed space-time formulation for the Dirichlet scalar wave problem,
but it also highlights that the Neumann scalar wave problem belongs to the same
family of problems in a dual fashion via the Hodge star operator. We believe
this insight is particularly meaningful when aiming to have space-time formulations
that allow for coupling of elliptic, parabolic and/or hyperbolic equations.

While we illustrated how these tools extend existing
results for the particular cases of Theorems~\ref{thm:ex_uni_L2} and
\ref{thm:no_uniqueness}, it goes without saying that the same mechanism works
for most results that begin from the ODE $u''+\mu u=f$, including
\cite{HLU25, Zan25}, to name a few. We note, however, that while
Theorems~\ref{thm:ex_uni_L2} and \ref{thm:no_uniqueness} do not actually
require $\Omega$ to be contractible, the topology of $\Omega$ may play a role
on other results of interest. For instance, the extension of the
resonance-aware data norms from \cite{STZ26} would require special treatment of
the Neumann zero mode when $k=n$, as well as possible harmonic modes in other
degrees when $\Omega$ is not contractible.

It is also worth mentioning that our framework leads to a different vectorial
formulation than that from \cite{Hau21, Hau24}, and could potentially offer
better numerical properties than those in \cite{Hau21,HaZ23}.

Last but not least, this unified framework also highlights a distinction
between time-dependent spatial differential forms and differential forms defined
directly over the space-time cylinder $Q$. Indeed, while for $k=0$ we have the
natural identification
\begin{equation*}
 L^2(0,T;L^2\Lambda^0(\Omega)) \cong L^2\Lambda^0(Q)= L^2(Q),
\end{equation*}
as used for the scalar wave equation; for $1\leq k\leq n$ the space of space-time
$k$-forms additionally contains components involving the temporal differential
$dt$. More precisely,
\begin{align*}
 L^2\Lambda^k(Q)\cong\;& L^2(0,T;L^2\Lambda^k(\Omega))\oplus\,dt\wedge L^2(0,T;
 L^2\Lambda^{k-1}(\Omega)).
\end{align*}
Nevertheless, for $k=n$, time-dependent spatial $n$-forms can also be identified
with space-time $(n+1)$-forms through the exterior product with $dt$. Hence,
in that case we can also identify
\begin{equation*}
 L^2(0,T;L^2\Lambda^n(\Omega)) \cong L^2\Lambda^{n+1}(Q)= L^2(Q).
\end{equation*}
Ongoing and future research deals with understanding how these different
representations influence space-time variational formulations for $0<k<n$.

\bibliographystyle{plain}
\bibliography{bib_HL}
 
\end{document}